\documentclass[12pt]{amsart}

\usepackage{amssymb,amscd}

\usepackage{enumerate}

\usepackage{graphicx}

\usepackage{hyperref}

\hypersetup{breaklinks=true}

\usepackage[all]{xy} 
\usepackage{mathtools} 

\usepackage{color} 

\makeatletter
\@namedef{subjclassname@2010}{%
  \textup{2010} Mathematics Subject Classification}
\makeatother

\newtheorem{thm}{Theorem}[section]

\newtheorem{lemma}[thm]{Lemma}

\newtheorem{propo}[thm]{Proposition}

\newtheorem*{funthmGl}{First Fundamental Theorem for $\GlV$}

\theoremstyle{definition}
\newtheorem{defin}[thm]{Definition}

\numberwithin{equation}{section}

\newcommand{\R}{\mathbb{R}} 
\newcommand{\Gl}{\mathrm{Gl}} 
\newcommand{\Tor}{\mathrm{Tor}} 
\newcommand{\Id}{I} 

\def\C{\mathcal{C}}

\def\RR{\mathbb{R}}
\def\R2n{\mathbb{R}^{2n}}

\def\Z2{\mathbb{Z}/2 \mathbb{Z}}

\def\T{\mathcal{T}}

\renewcommand{\d}{\ensuremath{\mathrm{d}}}
\def\dd{\ensuremath{\mathrm{d}}}

\def\qed{\hfill $\square$}

\def\Auts0{\protect \mathrm{Aut}(s_0)}

\def\Gl{{\rm GL}}
\def\GlV{{\rm GL}(V)}

\begin{document}

\baselineskip=17pt

\title{Uniqueness of the torsion-curvature pair} 

\author[R. Mart\'inez-Boh\'orquez]{Ra\'ul~Mart\'inez-Boh\'orquez}
\address{Departamento de Didáctica de las Ciencias Experimentales y Matem\'{a}ticas \\ Universidad de Extremadura \\ E-06071 Badajoz, Spain}
\email{raulmb@unex.es}

\author[J. Navarro]{Jos\'e ~Navarro}
\address{Departamento de Matem\'{a}ticas \\ Universidad de Extremadura \\ E-06071 Badajoz, Spain}
\email{navarrogarmendia@unex.es}

\author[J. B. Sancho]{Juan B. Sancho} 
\email{jsancho@unex.es}

\thanks{This project has been supported to the 85\% by the EU-ERDF, Junta de Extremadura, Autoridad de Gestión and Ministerio de Hacienda, through project GR24068. R.M.-B. and J.N. have also been supported by FEDER/AEI/MICINN through the grant PID2022-142024NB-I00 (SIGA). 
J.N. has also been supported by FEDER/AEI/MICINN through the grant PID2022-137667NA-I00 (GENTLE)}

\date{\today}

\begin{abstract}

On smooth manifolds of dimension $n \ge 4$, we prove that the torsion and curvature are, up to a scalar factor, the only pair of a vector-valued  2-form and an endomorphism-valued 2-form naturally associated with a linear connection that satisfy both the linear and differential Bianchi identities.

This result extends to arbitrary linear connections a recent characterisation of the curvature tensor of a symmetric linear connection obtained in (\cite{GMN_RACSAM}).
\end{abstract}

\maketitle


\section{Introduction}
\label{sec:intro}

Since the early days of the twentieth century, the study of natural operations has shed light on the role of fundamental constructions in Differential Geometry. One of the first examples is the characterization, due to E. Cartan \cite{CARTAN}, of the Einstein tensor of a pseudo-Riemannian metric, that provided a deeper understanding of the fundamental character of Einstein's equations of General Relativity (see also \cite{NS_JGP}).


Later on, this theory of natural constructions was developed so as to produce explicit descriptions of all natural operations of various kinds. As an example, P. B. Gilkey (\cite{GILKEY}) characterised the characteristic classes of Riemannian manifolds, a result that was later used to produce an alternative proof of the index theorem (\cite{ABP}). 
By the end of the last century, the monograph by Kol\'a\v{r}-Michor-Slov\'ak (\cite{KMSBOOK}) developed powerful mathematical methods and enhanced this theory to a greater level of generality.

In recent years, the theory has seen a resurgence of interest, fuelled by novel applications in physics (\cite{KMCMM}, \cite{KMM}) as well as by advances in different fields (\cite{BERNIG}, \cite{FHBAMS},  \cite{GMN_FEDOSOV}, \cite{NS_JGP}). Among these, we highlight a recent characterization of the curvature tensor of a symmetric, linear connection $\nabla$, as the only endomorphism-valued  natural 2-form $R$ that satisfies Bianchi's identities: 
$$ R \wedge I = 0 \quad , \quad \d _ \nabla R = 0 \ , $$
where $\Id : T X \to T X$ denotes the identity endomorphism (\cite[Thm. 3.13]{GMN_RACSAM}).


Our aim in this paper is to generalise this result to include general linear connections, which may have non-zero torsion.



In general, associated with any linear connection $\nabla$, two key tensors arise: its torsion tensor $\Tor$, which is a vector-valued 2-form, and its curvature tensor $R$, an endomorphism-valued 2-form. These tensors are intrinsically coupled via Bianchi's identities (see \cite[Thm. 2.5]{KN} for details): 
\begin{align}
    \d_\nabla \Tor &= R \wedge \Id \label{eq:bianchi1} \\
    \d_\nabla R &= 0 \label{eq:bianchi2} \ .
\end{align}
Our main result, Theorem~\ref{MAIN_THM}, establishes that these tensors, $\Tor$ and $R$, are the only pair of natural tensors fulfilling the system \eqref{eq:bianchi1}-\eqref{eq:bianchi2}, up to a common scalar factor. 

To prove this result, we  first have to describe the space of endomorphism-valued, natural 2-forms that are closed, which is interesting in its own right. We show that, when $n \ge 4$, this real vector space is 3-dimensional and it is spanned by the curvature tensor $R$ and two other tensors, constructed from traces of the curvature and the torsion (Theorem~\ref{TEOR_UNO}). 

Our proofs employ standard techniques from the theory of natural operations (\cite{GMN_RACSAM}, \cite{KMSBOOK}), 
which reduce the problem to describing vector spaces of $\Gl(n, \mathbb{R})$-equivariant linear maps between certain tensor spaces, defined using jets of the connection. On several occasions, we encounter the problem of checking that a number of tensors are linearly independent; these computations are routine, but too extensive to be detailed, so we use a computer program for these (cf. the proofs of Lemma \ref{L19}, Theorem \ref{TEOR_UNO} and Lemma \ref{LEMA_FIN}). 

\section{Naturality and normal tensors}\label{sec:statements} 

In this paper, $X$ will denote a smooth manifold of dimension $n$.

\subsection{Normal tensors}
\medskip

\begin{defin} 
A normal tensor of order $m\geq 0$ is a $\,(m+2,1)$-tensor $T$ on $X$ satisfying the
following symmetries:
\begin{itemize}
\item[-] it is symmetric in the last $m$ covariant indices; that is
\begin{equation}\label{S1}
T^l_{ijk_1\ldots k_m} = T^l_{ijk_{\sigma(1)} \ldots k_{\sigma(m)}} \ , \quad
\forall \ \sigma \in S_m\,;
\end{equation}

\item[-] the symmetrization of the $\,m+2\,$ covariant indices is zero:
\begin{equation}\label{S2}
\sum_{\sigma \in S_{m+2}} T^l_{\sigma (i) \sigma (j) \sigma (k_1)  \ldots \sigma (k_m)} = 0.
\end{equation} 
\end{itemize}
\end{defin} 

\begin{defin} Any linear connection $\nabla$ on $X$ has an associated sequence $\,\{N^m(\nabla)\}_{m\geq 0}\,$ of normal tensors, which are defined as follows: for any given point $x_0\in X$, consider normal coordinates $(x_1,\dots,x_n)$ with origin at $x_0$; the corresponding Christoffel symbols $\Gamma_{ij}^l$ satisfy the equation 
\begin{equation}\label{eq1}\sum x_ix_j\Gamma_{ij}^l\,=\,0\ .\end{equation}

The $m$-th normal tensor $N^m(\nabla)$ is defined at the point $x_0=0$ by the formula
$$N_{x_0}^m(\nabla)^l_{ijk_1\dots k_m}:=\,\frac{\partial^m\Gamma_{ij}^l}{\partial x_{k_1}\dots\partial x_{k_m}}(0)\ .$$

By successively partially differentiating in the equality (\ref{eq1}) and evaluating at $x_0=0$, it turns out that the tensors $N^m(\nabla)$ have the symmetries (\ref{S1}) and (\ref{S2}).
\end{defin}

These normal tensors $N^m(\nabla)$ are very useful for determining the tensors naturally associated with a connection (see Theorem \ref{NATURALES} below), but they are difficult to express in terms of covariant derivatives. For the purposes of this article, we only need the expressions for $N^0(\nabla)$ and $N^1(\nabla)$:

\begin{propo}\label{2.3} For any linear connection $\nabla$, it holds $\,N^0(\nabla)=\frac{1}{2}\mathrm{Tor}_\nabla$.
 \end{propo}
\begin{proof} Differentiating equation (\ref{eq1}) with $\partial x_i\partial x_j$  and evaluating at the origin produces
$$\Gamma_{ij}^l(0)+\Gamma_{ji}^l(0)\,=\,0\ .$$

Thus, the torsion tensor at the origin $x_0$ equals
$$\text{Tor}^l_{ij}\,=\, \Gamma_{ij}^l(0)-\Gamma_{ji}^l(0)\,=\,2\Gamma_{ij}^l(0)\,=\,2N_{ij}^l\ .$$
 
\qed \end{proof}

\begin{propo}\label{2.4}
The first normal tensor $N^1(\nabla)$ can be expressed in terms of the curvature tensor $R$, the torsion tensor $\mathrm{Tor}$ and its covariant derivative $\nabla\mathrm{Tor}$, using the following formula\footnote{Throughout the reminder of this paper, we will utilize Einstein's summation convention.}: 
 $$N_{ijk}^l= -\frac{1}{6}\left(-3R_{kij}^l+R_{jki}^l-R_{ijk}^l-2(\nabla \mathrm{Tor})_{ijk}^l-2(\nabla \mathrm{Tor})_{kji}^l+\mathrm{Tor}_{kj}^m\mathrm{Tor}_{mi}^l+\frac{1}{2}\mathrm{Tor}_{ij}^m\mathrm{Tor}_{km}^l \right)$$
 \end{propo}
\begin{proof} Let $(x_1,\dots,x_n)$ be normal coordinates with respect to a point $x_0$ and consider the corresponding Christoffel symbols $\Gamma_{ij}^l$. Recall the well-known expressions of $\mathrm{Tor}$, $\nabla \mathrm{Tor}$ and $R$ in terms of Christoffel symbols:
$$\begin{aligned}
\mathrm{Tor}_{ij}^k &=\Gamma_{ij}^k-\Gamma_{ji}^k \ , \\
(\nabla \mathrm{Tor})_{ijk}^l &=\mathrm{Tor}_{ij,k}^l-\Gamma_{ki}^m\mathrm{Tor}_{mj}^l-\Gamma_{kj}^m\mathrm{Tor}_{im}^l+\Gamma_{km}^l\mathrm{Tor}_{ij}^m \ , \\
R_{ijk}^l &=\Gamma_{jk,i}^l-\Gamma_{ik,j}^l+\Gamma_{jk}^m\Gamma_{im}^l-\Gamma_{ik}^m\Gamma_{jm}^l \ .
\end{aligned}$$

Now, let us consider the above formulas valued at the origin $x_0=0$, without denoting it explicitly. With this convention, we have by definition $N_{ij}^l=\Gamma_{ij}^l$ and $N_{ijk}^l=\Gamma_{ij,k}^l$. Furthermore, these tensors have the symmetries
$$\begin{aligned}
\Gamma_{ij}^l-\Gamma_{ji}^l=N_{ij}^k - N_{ji}^k = 0 \ , \\
N_{ijk}^l+N_{jki}^l+N_{kij}^l+N_{jik}^l+N_{ikj}^l+N_{kji}^l=0 \ .
\end{aligned}$$

Using these identities, the proof of the statement is routine.

\qed

\end{proof}

\subsection{Natural morphisms between tensor bundles}\phantom{.}\medskip

Let $\pi\colon T^{p,q}X=\bigotimes^pT^*X\otimes\bigotimes^qTX\longrightarrow X$ be the vector bundle of tensors of type $(p,q)$. Every diffeomorphism $\tau\colon U\to V$ between open sets of $X$ induces a diffeomorphism $\tau_* \colon T^{p,q}U\to T^{p,q}V$ that makes the following diagram commutative
$$\begin{CD}
T^{p,q}U @>{\tau_*}>{\sim}> T^{p,q}V\\
@V{\pi}VV @VV{\pi}V\\
U @>{\tau}>{\sim}> V
\end{CD}$$

\begin{defin} A morphism of vector bundles $\varphi\colon T^{p,q}X\to T^{\bar p, \bar q}X$ is  natural if, for every diffeomorphism $\tau\colon U\to V$ between open sets of $X$, the following  square is commutative
$$\begin{CD}
T^{p,q}U @>{\varphi}>> T^{\bar p, \bar q}U\\
@V{\tau_*}VV @VV{\tau_*}V\\
T^{p,q}V @>{\varphi}>> T^{\bar p, \bar q}V
\end{CD}$$

The vector space of natural morphisms of vector bundles $\varphi\colon T^{p,q}X\to T^{\bar p, \bar q}X$ will be denoted $\,\text{Hom}_{\text{nat}}(T^{p,q}X, T^{\bar p, \bar q}X)\,$.  
\end{defin}

This vector space is easy to compute: naturalness implies that every map $\varphi\colon T^{p,q}X\to T^{\bar p, \bar q}X$ is determined by its restriction to the fibre of a point $x_0\in X$ and, over such a fibre, the map is $\Gl (T_{x_0}X)$-equivariant; that is, there exists an $\mathbb{R}$-linear isomorphism
$$\begin{CD}
\text{Hom}_{\text{nat}}(T^{p,q}X, T^{\bar p, \bar q}X) @= \text{Hom}_{\Gl(T_{x_0}X)}(T^{p,q}_{x_0}X, T^{\bar p, \bar q}_{x_0}X) &\,=\, \text{Hom}_{ \Gl(T_{x_0}X)}(T^{p+\bar q,q+\bar p}_{x_0}X,\, \mathbb{R})\\ \varphi & \longmapsto &\varphi_{x_0} 
\end{CD}$$

The right-hand side of this equality can be described using the First Fundamental Theorem of invariant theory for the general linear group (\cite[Cor. 4.3.2]{GW} or \cite{COLLOQUIUM}):

\begin{funthmGl}\label{MainTheoremGl} Let $V$ be an $\mathbb{R}$-vector space of finite dimension, let $\GlV $ be the Lie group of its $\mathbb{R}$-linear automorphisms and let $T^{r,s}V=\bigotimes^rV^*\otimes\bigotimes^sV$.

The vector space $\,\mathrm{Hom}_{\GlV}\left( T^{r,s}V \, , \, \RR \right) \,$
 is zero if $r \neq s$, whereas, if $r = s$, it is spanned by the following ``total contractions":
$$ \phi_\sigma (\omega_1 \otimes \ldots\otimes\omega_r\otimes e_1\otimes \ldots \otimes e_r ) \, := \, \omega_1 (e_{\sigma (1)}) \cdot \ldots \cdot \omega_r (e_{\sigma (r)})\, , $$ 
where $\sigma$ runs over the permutations of $r$ elements.
\end{funthmGl}

As a direct consequence we obtain the following basic result:

\begin{thm}\label{2.6} 

 a) If $\,p-\bar p\neq q-\bar q\,$ then $\,\mathrm{Hom}_{\mathrm{nat}}(T^{p,q}X, T^{\bar p, \bar q}X)=0$.\smallskip

b) If $\,p-\bar p= q-\bar q\,$ then $\,\mathrm{Hom}_{\mathrm{nat}}(T^{p,q}X, T^{\bar p, \bar q}X)\,$ is generated by morphisms which are a composition of the following three types of operators:
\begin{enumerate}
    \item Contraction of $r$  pairs of indices:
    $$ C_{i_1,\dots,i_r}^{j_1\dots j_r}:\, T^{p,q}X\longrightarrow T^{(p-r,q-r)}X \ ; $$
    \item Products by the identity endomorphism:
    $$ \otimes^rI :\, T^{(p-r,q-r)}X\longrightarrow T^{\bar p, \bar q}X\quad ,\quad T\longmapsto T\otimes I\otimes\overset{r}{\dots}\otimes I \ ; $$
    \item Permutation of indices:
    $$ (\alpha,\beta) :\, T^{\bar p, \bar q}X\longrightarrow T^{\bar p, \bar q}X\quad ,\quad T_{i_1\dots i_{\bar p}}^{j_1\dots j_{\bar q}}\longmapsto T_{i_{\alpha(1)}\dots i_{\alpha(\bar p)}}^{j_{\beta(1)}\dots j_{\beta(\bar q)}} \ , $$
    where $\alpha $ and $\beta$ are permutations of $\bar p$ and $\bar q$ indices, respectively.
\end{enumerate}
\end{thm}

\smallskip
More generally, let $E\subseteq T^{p,q}X$ be a natural vector subbundle (naturalness here means that, for every diffeomorphism $\tau\colon U\to V$ between open sets of $X$, we have $\tau_*(E_{|U})=E_{|V}$). Since the restriction map $\,\text{Hom}_{\Gl(T_{x_0}X)}(T^{p,q}_{x_0}X, T^{\bar p, \bar q}_{x_0}X)\longrightarrow \text{Hom}_{\Gl(T_{x_0}X)}(E_{x_0}, T^{\bar p, \bar q}_{x_0}X)\,$ is surjective, then so it is the restriction map $\,\text{Hom}_{\text{nat}}(T^{p,q}X, T^{\bar p, \bar q}X)\longrightarrow \text{Hom}_{\text{nat}}(E, T^{\bar p, \bar q}X)\,$. Therefore, {\it the previous proposition remains valid if $\,\mathrm{Hom}_{\mathrm{nat}}(T^{p,q}X, T^{\bar p, \bar q}X)\,$ is replaced by $\,\mathrm{Hom}_{\mathrm{nat}}(E, T^{\bar p, \bar q}X)$.}\medskip

\subsection{Natural tensors associated to a linear connection}






Let $\C$ and $\T^{p,q}$ denote the sheaves of linear connections and that of $(p,q)$-tensors on $X$, respectively.

\begin{defin}A morphism of sheaves $\phi \colon \C \longrightarrow \T^{p,q}$ is regular if it is ``smooth", in the sense that it transforms any smooth family of linear connections into a family of tensors that is also smooth (see \cite[Section 2]{GMN_RACSAM} for details).
\end{defin}

In this paper, a natural tensor will always refer to a tensor naturally associated with linear connections, in the following precise sense:


\begin{defin}\label{naturaltensorconnections}
A $(p,q)$-tensor naturally associated with linear connections is a regular morphism of sheaves $\, T \colon \mathcal{C} \longrightarrow \mathcal{T}^{p,q} \,$ satisfying the following condition of naturalness:
$$ T ( \tau^* \nabla ) \, = \, \tau^* T(\nabla) \ , $$ for any diffeomorphism $\tau \colon U \to V$ between open sets of $X$, and for any linear connection $\nabla \in \C (V)$.
\end{defin}

The main examples we are interested in are the torsion and curvature of a linear connection. This objects define natural tensors:
$$\begin{aligned}
\Tor &: \C \longrightarrow \Omega^2_{TX}\qquad ,\qquad
  \nabla\mapsto\text{Tor}_\nabla \\   R &: \C \longrightarrow \Omega^2_{\mathrm{End}(TX)}\quad ,\quad \nabla\mapsto R_\nabla\ ,
  \end{aligned}$$
where $\Omega^2_{TX}$ and $\Omega^2_{\mathrm{End}(TX)}$ denote the sheaves of vector valued and endomorphism-valued $2$-forms on $X$, respectively.\medskip

Let us denote $\mathbf{N}^m\to X$ the fibre bundle of normal tensors of order $m$ over $X$. The proof of the following theorem, which is pivotal to our proofs, can be found in \cite[Thm. 2.3]{GMN_RACSAM}:

\begin{thm}\label{NATURALES}
There exists an $\mathbb{R}$-linear isomorphism
$$
\begin{CD}
\left\{
\begin{array}{c}
 \text{Natural tensors} \ \\
  \C \longrightarrow \mathcal{T}^{p,q} \ 
\end{array} \right\} @=
\bigoplus \limits_{d_0, \ldots , d_k} \mathrm{Hom}_{\mathrm{nat}}(S^{d_0}\mathbf{N}^0 \otimes \ldots \otimes S^{d_k}\mathbf{N}^k , T^{p,q} X ) \ ,
\end{CD}
$$
where $d_0, \ldots , d_k$ run over the non-negative integer solutions of the equation
\[
d_0 + \ldots + (k+1)d_k =p-q \ .
\] 

\end{thm}

In that theorem, a natural morphism of bundles $\varphi\colon S^{d_0}\mathbf{N}^0\otimes\cdots\otimes S^{d_k}\mathbf{N}^k\longrightarrow T^{p,q}X\,$ defines the following natural tensor $T\colon\mathcal{C}\to\mathcal{T}^{p,q}$, $$T(\nabla)\,=\varphi(N^0(\nabla),\overset{d_0}{\dots},N^0(\nabla),\dots,N^k(\nabla),\overset{d_k}{\dots},N^k(\nabla))\ .$$

\section{Main results}

As a first step, we will focus our attention on endomorphism-valued, natural 2-forms  that are closed. The curvature operator is an example of these, due to the so-called second, or differential, Bianchi identity (\ref{eq:bianchi2}).

On the other hand, 
if $\omega$ is an (ordinary) closed  $2$-form and $I: T X \to T X$ denotes the identity endomorphism, then $\omega \otimes I$ is an endomorphism-valued $2$-form that is also closed:
$$ \dd_\nabla ( \omega \otimes I ) = \d \omega \otimes I + \omega \otimes \dd_\nabla I = 0 + 0 = 0 \ .$$

As examples, both the first Chern form operator, $C_3^1 R=\text{tr}\circ R$ (or $R^k_{ijk}$ in coordinates), and the differential of the trace of the torsion, $\dd (C_1^1 \Tor)$, are closed $2$-forms (here, $C_i^j$ is the contraction operator of the $i$-th covariant index with the $j$-th contravariant index). Hence the corresponding endomorphism-valued forms $C_3^1 R \otimes I$ and $\dd (C_1^1 \Tor) \otimes I$ are also closed. 

We will prove that $R$, $C_3^1 R \otimes I$ and $\dd (C_1^1 \Tor) \otimes I$ are, essentially, the only natural endomorphism-valued 2-forms that are closed.

As a first step, let us compute the dimension of the vector space of endomorphism-valued, natural $2$-forms; that is to say, that of regular and natural morphisms of sheaves $  \C  \, \longrightarrow \ \Omega^2_{\mathrm{End}(TX)}$. 

From now on, we assume that $X$ is of dimension $n \geq 4$.

\begin{lemma}\label{L19}
The vector space of endomorphism-valued, natural $2$-forms has dimension 19.
\end{lemma}

\begin{proof} The sheaf $\Omega^2_{\mathrm{End}(TX)}$ is the sheaf of sections of the vector bundle $\Lambda^2T^*X\otimes T^*X\otimes TX$ which, in turn, is a natural vector subbundle of $ T^{3,1}X$. By Theorem \ref{NATURALES}, we have
$$
\begin{CD}
\left\{
\begin{array}{c}
 \text{Natural tensors} \ \\
  \C \longrightarrow \mathcal{T}^{3,1} \ 
\end{array} \right\} @=
\bigoplus \limits_{d_0, \ldots , d_k} \mathrm{Hom}_{\mathrm{nat}}(S^{d_0}\mathbf{N}^0 \otimes \ldots \otimes S^{d_k}\mathbf{N}^k ,\,  T^{3,1}X)\\
\bigcup & & \bigcup\\
\left\{
\begin{array}{c}
 \text{Natural tensors} \ \\
  \C \longrightarrow \Omega^2_{\mathrm{End}(TX)} \ 
\end{array} \right\} @=
\bigoplus \limits_{d_0, \ldots , d_k} \mathrm{Hom}_{\mathrm{nat}}(S^{d_0}\mathbf{N}^0 \otimes \ldots \otimes S^{d_k}\mathbf{N}^k , \,\Lambda^2T^*X\otimes T^*X\otimes TX) \ ,
\end{CD}
$$
where $d_0, \ldots , d_k$ run over the non-negative integer solutions of the equation
\[
d_0 + 2 d_1 + \ldots + (k+1)d_k =3-1=2 \ . 
\] 

There are only two solutions to the equation above, so let us analyze them separately:

\medskip 
\underline{Solution $d_1=1$, $ d_0 = d_2 = \ldots = 0$:} In this case we are led to compute generators for the vector space 
$$ \text{Hom}_{\text{nat}}(\mathbf{N}^1,\,\Lambda^2T^*X\otimes T^*X\otimes TX) \ . $$

To this end, we apply Theorem \ref{2.6}b: if $N_{ijk}^l$ is a  normal tensor in $\mathbf{N}^1$, we define four $(3,1)$-tensors $C_0, C_1, C_2$ and $C_3$ making index contractions and tensor products by the identity endomorphism:
\begin{itemize}
\item[-] $(C_0)_{ijk}^l:= N_{ijk}^l$
\item[-] $(C_1)_{ijk}^l:= N_{mij}^m \delta_k^l \ ,$ 
\item[-] $(C_2)_{ijk}^l:= N_{imj}^m \delta_k^l \ ,$ 
\item[-] $(C_3)_{ijk}^l:= N_{ijm}^m \delta_k^l \ .$ 
\end{itemize}

 Now, permutations of the covariant indices and skew-symmetrizations of indexes $i$ and $j$ in these four tensors produce a system of 12 generators:
\begin{itemize}
\item[-] $(T_{3\alpha+1})_{ijk}^l:=(C_\alpha)_{ijk}^l - (C_\alpha)_{jik}^l \ ,$
\item[-] $(T_{3\alpha+2})_{ijk}^l:=(C_\alpha)_{jki}^l - (C_\alpha)_{ikj}^l \ ,$
\item[-] $(T_{3\alpha+3})_{ijk}^l:=(C_\alpha)_{kij}^l - (C_\alpha)_{kji}^l \ ,$
\end{itemize}
where $\alpha=0,1,2,3$. 

However, we can remove the last item of the list, $T_{12}$, as it is a linear combination of $T_5, T_6, T_8, T_9$ and $T_{11}$, due to the symmetry (\ref{S2}) of the elements in $\mathbf{N}^1$.

\medskip
\underline{Solution $d_0=2$, $d_1 = \ldots = d_k = 0$:} In this case we have to analyze the vector space 
$$\text{Hom}_{\text{nat}}(S^2 \mathbf{N}^0,\,  \Lambda^2 T^* X \otimes T^* X \otimes T X) \ . $$

As in the previous case, given an element $N_{ij}^k N_{ab}^c$ of $\,S^2 \mathbf{N}^0$, we define the following five tensors, making index contractions and tensor products by the identity endomorphism (recall that $N_{ij}^l = - N_{ji}^l$):
\begin{itemize}
\item[-] $(D_1)_{ijk}^l:= N_{ij}^m N_{mk}^l \ ,$
\item[-] $(D_2)_{ijk}^l:= N_{ij}^l N_{mk}^m \ ,$
\item[-] $(D_{3})_{ijk}^l:= N_{si}^m N_{mk}^s \delta_j^l \ ,$
\item[-] $(D_{4})_{ijk}^l:= N_{ij}^m N_{ms}^s \delta_k^l \ ,$
\item[-] $(D_{5})_{ijk}^l:= N_{mi}^m N_{sj}^s \delta_k^l \ .$
\end{itemize}

Again, making permutations of covariant indices and skew-symmetrizations of the indices $i$ and $j$, we obtain the following system of eight generators for this second case:
\begin{itemize}
\item[-] $(T_{12})_{ijk}^l:=2(D_1)_{ijk}^l \ ,$
\item[-] $(T_{13})_{ijk}^l:=(D_1)_{jki}^l - (D_1)_{ikj}^l \ ,$
\item[-] $(T_{14})_{ijk}^l:=2(D_2)_{ijk}^l \ ,$
\item[-] $(T_{15})_{ijk}^l:=(D_2)_{jki}^l - (D_2)_{ikj}^l \ ,$
\item[-] $(T_{16})_{ijk}^l:=(D_3)_{jki}^l - (D_3)_{ikj}^l \ ,$
\item[-] $(T_{17})_{ijk}^l:=2(D_4)_{ijk}^l \ ,$
\item[-] $(T_{18})_{ijk}^l:=(D_4)_{jki}^l - (D_4)_{ikj}^l \ ,$
\item[-] $(T_{19})_{ijk}^l:=2(D_5)_{ijk}^l \ .$
\end{itemize}

\medskip 
Summing up, we have a total of 19 generators, and it remains to prove that they are linearly independent. It suffices to check this on a particular affine manifold, and using the following example will be enough.

\medskip \underline{Test example:} Let $\nabla$ be the affine connection on $\RR^4$ defined, in cartesian coordinates $x_1,x_2,x_3,x_4$, by the following Christoffel symbols\footnote{Any Christoffel symbol that is not mentioned is assumed to be null.}:
\begin{itemize}
\item[-] $\Gamma_{12}^1=x_3\ ,$
\item[-] $\Gamma_{43}^3=x_1 x_4\ ,$
\item[-] $\Gamma_{31}^3=x_2 x_4\ .$
\end{itemize}  

The tensors $N^0(\nabla)$ and $N^1(\nabla)$ are computed using propositions \ref{2.3} and \ref{2.4}. Now, a routine calculation, that was carried out using a computer program, shows that those 19 tensors above, $T_1, \ldots , T_{19}$, are $\RR$-linearly independent.

\qed \end{proof}

\begin{thm}\label{TEOR_UNO}
Let $X$ be a smooth manifold of dimension $n\geq 4$. The vector space of endomorphism-valued, closed $2$-forms naturally associated with linear connections, has dimension $3$, and it is spanned by $R$, $C_3^1 R \otimes I$ and $\dd (C_1^1 \Tor) \otimes I$.
\end{thm}

\begin{proof} We already mentioned at the beginning of the section that $R$, $C_3^1 R \otimes I$ and $\dd (C_1^1 \Tor) \otimes I$ are closed 2-forms. Let us see that there are no other closed, natural 2-forms that are linearly independent of those three.

Let us consider the generators $T_1, \ldots , T_{19}$ obtained in Lemma \ref{L19}. With the aid of a computer program, we differentiate these nineteen 2-forms in the aforementioned test example, and find exactly 
three linear combinations whose differentials are null, which are
$$ T_{2} -  T_{13} \quad , \quad T_{4} \quad , \quad  T_{7} \ .$$

These linear combinations correspond, respectively, to the natural 2-forms:
$$ R \quad , \quad - \dd (C_1^1 \Tor) \otimes I - C_3^1 R \otimes I \quad , \quad - C_3^1 R \otimes I \ . $$ \qed
\end{proof}

In the following two lemmata, the identity endomorphism $I$ is considered as a vector-valued 1-form. The next result can be can be found in \cite[Lemma 3.5]{GMN_RACSAM}. Recall that we assume that the dimension of $X$ is greater or equal than four.


\begin{lemma}\label{LEMA_FIN}
The torsion tensor $\Tor_\nabla$, and the 2-form $H_\nabla := C_1^1 (\Tor_\nabla ) \wedge I $, are a basis for the vector space of vector-valued natural 2-forms.
\end{lemma}

\begin{lemma}\label{3.4}
The four natural vector-valued 3-forms $R \wedge I$, $\,(C_3^1 R \otimes I) \wedge I$, $\,(\dd (C_1^1 \Tor) \otimes I) \wedge I\,$ and $\,\mathrm{d}_\nabla H\,$ are $\mathbb{R}$-linearly independent.\end{lemma}

\begin{proof} Again, it is enough to check it in the test example introduced in the proof of Lemma \ref{L19}, and we implemented this routine computation in a computer program.
\qed\end{proof}

\begin{thm}\label{MAIN_THM}
Let $X$ be a smooth manifold of dimension $n\geq 4$ and let $\alpha , \beta $ be a pair of a vector-valued 2-form and an endomorphism-valued 2-form, naturally associated with a linear connection $\nabla$. If these forms $\alpha$ and $ \beta$ satisfy the Bianchi identities 
$$\dd_\nabla \alpha = \beta \wedge I \qquad \mbox{ and } \qquad \dd_\nabla \beta = 0 \ , $$
then $ \alpha = \lambda\, \Tor_\nabla$ and $\beta = \lambda\, R_\nabla$, for some $\lambda \in \RR$.
\end{thm}
\begin{proof} On the one hand, by Theorem \ref{TEOR_UNO} there exist real numbers $\lambda_1 , \lambda_2$ and $\lambda_3$ such that
$$\beta\,=\,\lambda_1 R+\lambda_2\,C_3^1R\otimes I+\lambda_3\,\text{d}(C_1^1\text{Tor})\otimes I\ ,$$
and hence
$$\beta\wedge I\,=\,\lambda_1 R\wedge I +\lambda_2\,(C_3^1R\otimes I)\wedge I+\lambda_3\,(\text{d}(C_1^1\text{Tor})\otimes I)\wedge I\ .$$

On the other hand, due to Lemma \ref{LEMA_FIN}, there exist  $\lambda , \mu \in \RR$ such that 
$$ \alpha = \lambda \Tor + \mu H \ ,$$ 
and, differentiating, we also obtain
$$\beta\wedge I\,=\,\text{d}_\nabla\alpha\,=\,\lambda\text{d}_\nabla\text{Tor}+\mu\text{d}_\nabla H\,=\, \lambda R\wedge I+\mu\text{d}_\nabla H\ .$$

Upon comparing the two resulting expressions for $\beta\wedge I$, and using Lemma \ref{3.4}, we find
$$\lambda_1\,=\,\lambda\quad ,\quad \lambda_2\,=\,\lambda_3\,=\,\mu\,=\, 0\ ,$$
so that $\alpha=\lambda\text{Tor}\,$ and $\,\beta=\lambda R$.

\qed

\end{proof}

A related result is obtained in (\cite{NS_DGA} Thm. A2 and Cor. A5), where the curvature 2-form of a connection on a principal bundle is characterised in terms of a different notion of naturalness.

\noindent \textbf{Conflicts of interest:} The authors declare no conflicts of interest.

\section*{Appendix: MATLAB script}

The MATLAB script utilised for the computations can be accessed at the following GitHub repository: 

https://github.com/RaulMartinezBohorquez/Uniqueness-of-the-Torsion-Curvature-Pair-of-Operators

\end{document}